# Generalization of Ehrlich-Kjurkchiev method for multiple roots of algebraic equations

ANTON I. ILIEV

In this paper a new method which is a generalization of the Ehrlich-Kjurkchiev method is developed. The method allows to find simultaneously all roots of the algebraic equation in the case when the roots are supposed to be multiple with known multiplicities. The offered generalization does not demand calculation of derivatives of order higher than first simultaneously keeping quaternary rate of convergence which makes this method suitable for application from the practical point of view.

**Introduction.** The problem of simultaneous finding all roots (SFAR) of polynomial equation is connected with the works of Weierstrass [1] and Dochev [2]. Their formula represents the modification of the Newton method for individual determination of the roots of the algebraic equation. The method suggested by Weierstrass-Dochev has quadratic rate of convergence. Using divided differences with multiple knots Semerdzhiev [3] generalizes the method of Weierstrass-Dochev, when the roots are multiple. Another approach to SFAR is proposed by Ehrlich [4]. He shows that for his method the modification of Gauss-Seidel is applicable. The Ehrlich method is generalized for the case of multiple roots of algebraic, trigonometric and exponential equations in [5]. The presented there method is of the same complexity as the method of Ehrlich for simple roots. The same method is also generalized [6] for the case of polynomial equations on arbitrary Chebyshev system having multiple roots with known multiplicities. The rate of convergence is also cubic, but at each step determinants should be calculated.

In this paper a new method representing a generalization of the method of Ehrlich-Kjurkchiev [7] for algebraic equations is obtained. Our method has also quaternary rate of convergence and is effective from a computational point of view and it can be used for SFAR, when they have known multiplicities. This new method is of the same complexity as the method of Ehrlich-Kjurkchiev for simple roots.

Let the algebraic polynomial

(1) $\qquad A_n(x) = x^n + a_1 x^{n-1} + \ldots + a_n$

be given and $x_1, x_2, \ldots, x_m$ be its roots with known multiplicities $\alpha_1, \alpha_2, \ldots, \alpha_m$ $(\alpha_1 + \alpha_2 + \ldots + \alpha_m = n)$. When the roots are simple $(\alpha_1 = \alpha_2 = \ldots \alpha_m = 1)$ Kjurkchiev [7] proposed the following iteration formula for SFAR

$$x_i^{[k+1]} = x_i^{[k]} - A_n(x_i^{[k]}) / \left[ A_n'(x_i^{[k]}) - A_n(x_i^{[k]}) W_i'^{[k]}(x_i^{[k]}) / W_i^{[k]}(x_i^{[k]}) \right.$$

(2)
$$\left. + A_n(x_i^{[k]}) \sum_{j=1, j \neq i}^{n} A_n(x_j^{[k]}) (x_i^{[k]} - x_j^{[k]})^{-2} / W_j^{[k]}(x_j^{[k]}) \right]$$

$i = \overline{1,n}, k = 0,1,2,...$

where $W_p^{[k]}(x_p^{[k]}) = \prod_{l=1, l \neq p}^{n} (x_p^{[k]} - x_l^{[k]})$, $p = \overline{1,n}$.

Our principal aim is to generalize (2) for the case when the roots are not simple, i.e. when $\alpha_1, \alpha_2, ..., \alpha_m$ are arbitrary.

So, we define

(3) $$x_i^{[k+1]} = x_i^{[k]} - \alpha_i \left[ S_i^{[k]}(x_i^{[k]}) + \sum_{j=1, j \neq i}^{m} \alpha_j (x_j^{[k]} - x_i^{[k]})^{-2} A_n(x_j^{[k]}) \left[ S_j^{[k]}(x_j^{[k]}) / \alpha_j \right]^{\alpha_j - 1} / Q_j^{[k]}(x_j^{[k]}) \right]^{-1}$$

$i = \overline{1,m}, k = 0,1,2,...$

where $S_p^{[k]}(x_p^{[k]}) = A_n'(x_p^{[k]}) / A_n(x_p^{[k]}) - Q_p'^{[k]}(x_p^{[k]}) / Q_p^{[k]}(x_p^{[k]})$, $p = \overline{1,m}$

and $Q_p^{[k]}(x_p^{[k]}) = \prod_{l=1, l \neq p}^{m} (x_p^{[k]} - x_l^{[k]})^{\alpha_l}$, $p = \overline{1,m}$.

The following theorem shows, that the iterative formula (3) is convergent to the roots with quaternary rate of convergence.

**Theorem.** *Let $c$, $q$ and $d \stackrel{\text{def}}{=} \min_{i \neq j} |x_i - x_j|$ be real positive constants, so that the following inequalities be fulfilled*

$q < 1$, $d - 2c > 0$, $2c^2 n(d-2c)^{-2} \left[ c(d-2c)^{-1} + (1 + c(d-2c)^{-1})[N + MN + M] \right] < \alpha_i$, $i = \overline{1,m}$

*where* $M \stackrel{\text{def}}{=} \left[ [1 + c/(d-2c)]^n - 1 \right]$ *and* $N \stackrel{\text{def}}{=} \left[ [1 + n(c/(d-2c))^2]^{n-1} - 1 \right]$. *It is assumed that the initial approximations $x_1^{[0]}, x_2^{[0]}, ..., x_m^{[0]}$ to the roots of* (1) $x_1, x_2, ..., x_m$ *are chosen so that the inequalities $|x_i^{[0]} - x_i| < cq$, $i = \overline{1,m}$ hold true. Then for every natural number $k$ the inequalities*

(4) $|x_i^{[k]} - x_i| < cq^{4^k}$, $i = \overline{1,m}$

*are satisfied.*

**Proof.** Let us assume that the inequalities (4) be fulfilled for some natural $k$. It will be proved, that $|x_i^{[k+1]} - x_i| < cq^{4^{k+1}}$, $i = \overline{1,m}$.

Easily one can find that

(5) $Q_i'^{[k]}(x_i^{[k]}) / Q_i^{[k]}(x_i^{[k]}) = \sum_{j=1, j \neq i}^{m} \alpha_j / (x_i^{[k]} - x_j^{[k]})$, $i = \overline{1,m}$

and

(6) $S_i^{[k]}(x_i^{[k]}) = \alpha_i / (x_i^{[k]} - x_i) + \sum_{j=1, j \neq i}^{m} \alpha_j (x_j - x_j^{[k]}) [(x_i^{[k]} - x_j^{[k]})(x_i^{[k]} - x_j)]^{-1}, i = \overline{1, m}.$

If we subtract $x_i$ from both sides of the iterative formula (3), then using (5),(6) and reducing under common denominator the right side of (3) we receive

(7) $x_i^{[k+1]} - x_i = \left[ \alpha_i + (x_i^{[k]} - x_i) \sum_{j=1, j \neq i}^{m} \alpha_j (x_j - x_j^{[k]}) [(x_i^{[k]} - x_j^{[k]})(x_i^{[k]} - x_j)]^{-1} \right.$

$\left. + (x_i^{[k]} - x_i) P_i^{[k]}(x_i^{[k]}) - \alpha_i \right] [S_i^{[k]}(x_i^{[k]}) + P_i^{[k]}(x_i^{[k]})]^{-1}, i = \overline{1, m}$

where $P_i^{[k]}(x_i^{[k]}) = \sum_{j=1, j \neq i}^{m} \alpha_j A_n(x_j^{[k]})(x_j^{[k]} - x_i^{[k]})^{-2} [S_j^{[k]}(x_j^{[k]}) / \alpha_j]^{\alpha_j - 1} / Q_j^{[k]}(x_j^{[k]}).$

Further, if from the numerator and the denominator of (7) factors $(x_i^{[k]} - x_i)$ and $(x_i^{[k]} - x_i)^{-1}$ respectively are separated, then (7) can be written in the form

(8) $x_i^{[k+1]} - x_i = (x_i^{[k]} - x_i)^2 \left[ \sum_{j=1, j \neq i}^{m} \alpha_j (x_j - x_j^{[k]}) [(x_i^{[k]} - x_j^{[k]})(x_i^{[k]} - x_j)]^{-1} + P_i^{[k]}(x_i^{[k]}) \right]$

$\times \left[ \alpha_i + (x_i^{[k]} - x_i) \sum_{j=1, j \neq i}^{m} \alpha_j (x_j - x_j^{[k]}) [(x_i^{[k]} - x_j^{[k]})(x_i^{[k]} - x_j)]^{-1} + (x_i^{[k]} - x_i) P_i^{[k]}(x_i^{[k]}) \right]^{-1}$

$i = \overline{1, m}.$

Now we transform $P_i^{[k]}(x_i^{[k]}), i = \overline{1, m}$ by the following way

$P_i^{[k]}(x_i^{[k]}) = \sum_{j=1, j \neq i}^{m} \alpha_j (x_j^{[k]} - x_i^{[k]})^{-2} [S_j^{[k]}(x_j^{[k]}) / \alpha_j]^{\alpha_j - 1} \prod_{l=1}^{m} (x_j^{[k]} - x_l)^{\alpha_l} \prod_{s=1, s \neq j}^{m} (x_j^{[k]} - x_s^{[k]})^{-\alpha_s}$

$= \sum_{j=1, j \neq i}^{m} \alpha_j (x_j^{[k]} - x_j)(x_i^{[k]} - x_j^{[k]})^{-2} \left[ 1 + ((x_j^{[k]} - x_j) / \alpha_j) \sum_{l=1, l \neq j}^{m} \alpha_l (x_l - x_l^{[k]}) [(x_j^{[k]} - x_l)(x_j^{[k]} - x_l^{[k]})]^{-1} \right]^{\alpha_j - 1}$

$\times \prod_{s=1, s \neq j}^{m} (x_j^{[k]} - x_s^{[k]})^{-\alpha_s} (x_j^{[k]} - x_s)^{\alpha_s}, i = \overline{1, m}.$

If we multiply the numerator and the denominator of the previous equation with $(x_i^{[k]} - x_j)$ and put the obtained expression for $P_i^{[k]}(x_i^{[k]}), i = \overline{1, m}$ into (8) we receive

$$x_i^{[k+1]} - x_i = \left(x_i^{[k]} - x_i\right)^2 \sum_{j=1, j \neq i}^{m} \alpha_j \left(x_j - x_j^{[k]}\right) \left[\left(x_i^{[k]} - x_j^{[k]}\right)\left(x_i^{[k]} - x_j\right)\right]^{-1}$$

$$\times \left[1 - \left(x_i^{[k]} - x_j\right)\left(x_i^{[k]} - x_j^{[k]}\right)^{-1} R_j^{[k]}\left(x_j^{[k]}\right) \prod_{s=1, s \neq j}^{m} \left(x_j^{[k]} - x_s^{[k]}\right)^{-\alpha_s} \left(x_j^{[k]} - x_s\right)^{\alpha_s}\right]$$

(9)
$$\times \left[\alpha_i + \left(x_i^{[k]} - x_i\right) \sum_{j=1, j \neq i}^{m} \alpha_j \left(x_j - x_j^{[k]}\right)\left[\left(x_i^{[k]} - x_j^{[k]}\right)\left(x_i^{[k]} - x_j\right)\right]^{-1} \left[1 - \left(x_i^{[k]} - x_j\right)\left(x_i^{[k]} - x_j^{[k]}\right) R_j^{[k]}\left(x_j^{[k]}\right)\right.\right.$$

$$\left.\left.\times \prod_{s=1, s \neq j}^{m} \left(x_j^{[k]} - x_s^{[k]}\right)^{-\alpha_s} \left(x_j^{[k]} - x_s\right)^{\alpha_s}\right]\right]^{-1}, i = \overline{1, m}$$

where $R_j^{[k]}\left(x_j^{[k]}\right) = \left[1 + \left(\left(x_j^{[k]} - x_j\right) / \alpha_j\right) \sum_{l=1, l \neq j}^{m} \alpha_l \left(x_l - x_l^{[k]}\right)\left[\left(x_j^{[k]} - x_l\right)\left(x_j^{[k]} - x_l^{[k]}\right)\right]^{-1}\right]^{\alpha_j - 1}$.

Denote

(10) $Y_{ij}^{[k]}\left(x_i^{[k]}, x_j^{[k]}\right) \overset{\text{def}}{=} 1 - \left(x_i^{[k]} - x_j\right)\left(x_i^{[k]} - x_j^{[k]}\right)^{-1} R_j^{[k]}\left(x_j^{[k]}\right) \prod_{s=1, s \neq j}^{m} \left(x_j^{[k]} - x_s^{[k]}\right)^{-\alpha_s} \left(x_j^{[k]} - x_s\right)^{\alpha_s}$

$i, j = \overline{1, m}, i \neq j$.

The expression for $R_j^{[k]}\left(x_j^{[k]}\right), j = \overline{1, m}$ can be transformed in the form

(11)
$$R_j^{[k]}\left(x_j^{[k]}\right) = 1 + \sum_{r=1}^{\alpha_j - 1} (\alpha_j - 1)! \left[r!(\alpha_j - 1 - r)!\right]^{-1}$$

$$\times \left[\left(\left(x_j^{[k]} - x_j\right)/\alpha_j\right) \sum_{l=1, l \neq j}^{m} \alpha_l \left(x_l - x_l^{[k]}\right)\left[\left(x_j^{[k]} - x_l\right)\left(x_j^{[k]} - x_l^{[k]}\right)\right]^{-1}\right]^r, j = \overline{1, m}.$$

Using the following inequalities

(12)
$$\left|x_j^{[k]} - x_l\right| \geq \left|x_j - x_l\right| - \left|x_j - x_j^{[k]}\right| \geq d - cq^{4^k} > d - c > d - 2c$$

$$\left|x_j^{[k]} - x_l^{[k]}\right| \geq \left|x_j^{[k]} - x_l\right| - \left|x_l - x_l^{[k]}\right| \geq d - 2cq^{4^k} > d - 2c, l, j = \overline{1, m}, l \neq j$$

the estimate for $\left|R_j^{[k]}\left(x_j^{[k]}\right)\right|, j = \overline{1, m}$, can be found.

Namely, from (11), (12) and $N \overset{\text{def}}{=} \left[\left[1 + n(c/(d - 2c))^2\right]^{n-1} - 1\right]$ it follows that

(13) $\left|R_j^{[k]}\left(x_j^{[k]}\right)\right| \leq 1 + \left(q^{4^k}\right)^2 N, j = \overline{1, m}$.

Further, we make the following transformation

$$Z_j^{[k]}\left(x_j^{[k]}\right) \stackrel{\text{def}}{=} \prod_{s=1, s \neq j}^{m} \left(x_j^{[k]} - x_s^{[k]}\right)^{-\alpha_s} \left(x_j^{[k]} - x_s\right)^{\alpha_s}$$

$$= \left[ \prod_{s=1, s \neq j}^{m} \left(x_j^{[k]} - x_s^{[k]}\right)^{-\alpha_s} \left(x_j^{[k]} - x_s\right)^{\alpha_s} - \prod_{s=2, s \neq j}^{m} \left(x_j^{[k]} - x_s^{[k]}\right)^{-\alpha_s} \left(x_j^{[k]} - x_s\right)^{\alpha_s} \right]$$

$$+ \ldots + \left[ \prod_{s=m-1, s \neq j}^{m} \left(x_j^{[k]} - x_s^{[k]}\right)^{-\alpha_s} \left(x_j^{[k]} - x_s\right)^{\alpha_s} - \prod_{s=m, s \neq j}^{m} \left(x_j^{[k]} - x_s^{[k]}\right)^{-\alpha_s} \left(x_j^{[k]} - x_s\right)^{\alpha_s} \right]$$

$$+ \left[ \prod_{s=m, s \neq j}^{m} \left(x_j^{[k]} - x_s^{[k]}\right)^{-\alpha_s} \left(x_j^{[k]} - x_s\right)^{\alpha_s} - 1 \right] + 1, \; j = \overline{1, m}.$$

Obviously, the expression $Z_j^{[k]}\left(x_j^{[k]}\right), j = \overline{1,m}$ can be presented in the form

$$Z_j^{[k]}\left(x_j^{[k]}\right) = 1 + \sum_{l=1, l \neq j}^{m} \left(x_j^{[k]} - x_l\right)\left(x_j^{[k]} - x_l^{[k]}\right)^{-\alpha_l} \prod_{s=l+1, s \neq j}^{m} \left(x_j^{[k]} - x_s^{[k]}\right)^{-\alpha_s} \left(x_j^{[k]} - x_s\right)^{\alpha_s}$$

(14)

$$\times \sum_{p=1}^{\alpha_l} \left(x_j^{[k]} - x_l\right)^{\alpha_l - p} \left(x_j^{[k]} - x_l^{[k]}\right)^{p-1}, \; j = \overline{1,m}.$$

On the other hand we have

(15) $Z_j^{[k]}\left(x_j^{[k]}\right) = \prod_{s=1, s \neq j}^{m} \left[1 + \left(x_s^{[k]} - x_s\right) / \left(x_j^{[k]} - x_s^{[k]}\right)\right]^{\alpha_s}, \; j = \overline{1,m}.$

Now, using (12), (15) and $M \stackrel{\text{def}}{=} \left[\left[1 + c/(d - 2c)\right]^n - 1\right]$ one can find the following estimate

(16) $\left|Z_j^{[k]}\left(x_j^{[k]}\right)\right| \leq 1 + q^{4^k} M, \; j = \overline{1,m}.$

In order to receive an estimate for $Y_{ij}^{[k]}\left(x_i^{[k]}, x_j^{[k]}\right), i, j = \overline{1,m}, i \neq j$ we transform this expression by the use of (11) and (14) as follows

$$Y_{ij}^{[k]}\left(x_i^{[k]}, x_j^{[k]}\right) = 1 - \left(x_i^{[k]} - x_j\right)\left(x_i^{[k]} - x_j^{[k]}\right)^{-1} - \left(x_i^{[k]} - x_j\right)\left(x_i^{[k]} - x_j^{[k]}\right)^{-1}$$

(17)

$$\times \left[ \sum_{r=1}^{\alpha_j - 1} (\alpha_j - 1)! \left[r!(\alpha_j - 1 - r)!\right]^{-1} \left[\left(\left(x_j^{[k]} - x_j\right)/\alpha_j\right) \sum_{l=1, l \neq j}^{m} \alpha_l\left(x_l - x_l^{[k]}\right)\left[\left(x_j^{[k]} - x_l\right)\left(x_j^{[k]} - x_l^{[k]}\right)\right]^{-1}\right]^{r}\right]$$

$$+ \sum_{l=1, l \neq j}^{m} \left(x_l^{[k]} - x_l\right)\left(x_j^{[k]} - x_l^{[k]}\right)^{-\alpha_l} \prod_{s=l+1, s \neq j}^{m} \left(x_j^{[k]} - x_s^{[k]}\right)^{-\alpha_s} \left(x_j^{[k]} - x_s\right)^{\alpha_s} \sum_{p=1}^{\alpha_l} \left(x_j^{[k]} - x_l\right)^{\alpha_l - p} \left(x_j^{[k]} - x_l^{[k]}\right)^{p-1}$$

$$\times \left[1 + \sum_{r=1}^{\alpha_j - 1} (\alpha_j - 1)! \left[r!(\alpha_j - 1 - r)!\right]^{-1} \left[\left(\left(x_j^{[k]} - x_j\right)/\alpha_j\right) \sum_{l=1, l \neq j}^{m} \alpha_l\left(x_l - x_l^{[k]}\right)\left[\left(x_j^{[k]} - x_l\right)\left(x_j^{[k]} - x_l^{[k]}\right)\right]^{-1}\right]^{r}\right]$$

$i, j = \overline{1, m}, i \neq j.$

Using (12), (13), (16) and (17) it is easy to receive the estimate

(18) $$\left|Y_{ij}^{[k]}\left(x_i^{[k]},x_j^{[k]}\right)\right|\le q^{4^k}\left[c(d-2c)^{-1}+\left(1+cq^{4^k}(d-2c)^{-1}\right)\left[q^{4^k}N+M\left[1+\left(q^{4^k}\right)^2N\right]\right]\right]$$

$i,j=\overline{1,m}, i\neq j.$

Now, from (9) we obtain

(19) $$\left|x_i^{[k+1]}-x_i\right|\le \left(x_i^{[k]}-x_i\right)^2\sum_{j=1,j\neq i}^{m}\alpha_j\left|x_j-x_j^{[k]}\right|\left[\left|x_i^{[k]}-x_j^{[k]}\right|\left|x_i^{[k]}-x_j\right|\right]^{-1}\left|Y_{ij}^{[k]}\left(x_i^{[k]},x_j^{[k]}\right)\right|$$

$$\times\left[\alpha_i-\left|x_i^{[k]}-x_i\right|\sum_{j=1,j\neq i}^{m}\alpha_j\left|x_j-x_j^{[k]}\right|\left[\left|x_i^{[k]}-x_j^{[k]}\right|\left|x_i^{[k]}-x_j\right|\right]^{-1}\left|Y_{ij}^{[k]}\left(x_i^{[k]},x_j^{[k]}\right)\right|\right]^{-1}, i=\overline{1,m}.$$

Finally, with the help of inequalities (12) the estimate (18) and using the conditions of the theorem from (19) we receive

$$\left|x_i^{[k+1]}-x_i\right|<c^2\left(q^{4^k}\right)^2 ncq^{4^k}(d-2c)^{-2}q^{4^k}\left[c/(d-2c)+\left(1+cq^{4^k}(d-2c)^{-1}\right)\left[q^{4^k}N+M\left[1+\left(q^{4^k}\right)^2N\right]\right]\right]$$

$$\times\left[\alpha_i-c^2\left(q^{4^k}\right)^2 n(d-2c)^{-2}q^{4^k}\left[c(d-2c)^{-1}+\left(1+cq^{4^k}(d-2c)^{-1}\right)\left[q^{4^k}N+M\left[1+\left(q^{4^k}\right)^2N\right]\right]\right]\right]^{-1}$$

$$<cq^{4^{k+1}}c^2n(d-2c)^{-2}\left[c/(d-2c)+\left(1+c(d-2c)^{-1}\right)[N+MN+M]\right]$$

$$\times\left[\alpha_i-c^2n(d-2c)^{-2}\left[c(d-2c)^{-1}+\left(1+c(d-2c)^{-1}\right)[N+MN+M]\right]\right]^{-1}<cq^{4^{k+1}}, i=\overline{1,m}.$$

Thus the theorem is completely proved.

**Remark.** In the case when $\alpha_1=\alpha_2=\ldots=\alpha_m=1$ the method (3) coincides with the method (2).

**Numerical example.** For the equation $A_6(x)=(x+2)^2(x-1)(x-3)^3=0$ using initial approximation $x_1^{[0]}=-3, x_2^{[0]}=0.1, x_3^{[0]}=4$ with the help of the method (3) the roots with accuracy of 18 decimal digits have been found only after 3 iterations.

| k | $x_1^{[k]}$ | $x_2^{[k]}$ | $x_3^{[k]}$ |
|---|---|---|---|
| 0 | -3.000000000000000000 | 0.100000000000000000 | 4.000000000000000000 |
| 1 | -1.989380609181193540 | 0.995064651338749428 | 3.026047103321694120 |
| 2 | -1.999999999677379630 | 0.999999994237752166 | 3.000000006833252880 |
| 3 | -2.000000000000000000 | 1.000000000000000000 | 3.000000000000000000 |

**Acknowledgements.** The author is deeply grateful to Dr. Khristo Semerdzhiev for many useful discussions concerning this paper.

*University of Plovdiv*
http://www.pu.acad.bg
*Faculty of Mathematics and Informatics*
http://www.fmi.pu.acad.bg
*Department of Numerical Methods*
*24 Tzar Assen str.*
*4000 Plovdiv, BULGARIA*
*e-mail:* aii@pu.acad.bg
*URL:* http://anton.iliev.tripod.com